\documentclass[11pt]{article}
\usepackage{graphicx}
\usepackage{url}
\begin{document}

\title{Matching sequences of two digits in matrices is hard}
\author{Nicolaos Matsakis\\
University of California- Irvine\\
nmatsaki@uci.edu}
\date{}
\maketitle

\begin{abstract}
We introduce a new -as far as we know- problem, according to which we are asked to find whether a given sequence of '0' and '1' digits is matched completely by a simple path in matrices having entries among those two digits (and others too) and prove that this problem is NP-complete. 
\end{abstract}

\section{Introduction}
We introduce the following problem: We are given a nxm matrix of 0's, 1's and any other digits and a sequence of 0's and 1's only. The length of the sequence is less or equal to the sum of the numbers of 0's and 1's in the provided matrix. We are asking whether there is a path in the matrix, when moving in horizontal, vertical or diagonal way and always between adjacent matrix entries, such that it is matched completely by the specified sequence, visiting entries of 0 and 1, no more than once each. As an example we have the matrix:

\[ \left( \begin{array}{cccccc}
1 & 0 & 1 & 1 & 0 &1\\
0 & 1 & 1 & 1 & 2 &0\\
0 & 0 & 0 & 1 & 1 &1 \\
1 & 0 & 3 & 0 & 0 &0 \\
1 & 1 & 1 & 0 & 3 &1\\
2 & 5 & 1 & 0 & 1 &1\\
2 & 4 & 0 & 0 & 0 & 0
\end{array} \right)\]    

\noindent and the sequence 1011010111110010010111000011101000, for which, apparently, there is a simple path in the matrix. A more restricted case of this problem can be retrieved from \cite{eppstein} (matrix is square and must have only 0's and 1's, while the sequence has length exactly equal to the number of entries). We may find various NP-complete problems regarding matrices of digits \cite{garey,pap}, such as the Consecutive Ones Submatrix problem and the Consequtive Block Minimization problem, however we cannot find this specific problem in a relevant source; so we name it 0-1 String to Matrix problem.

\section{Proof of NP-completeness}

\begin{figure}
\begin{center}
\includegraphics[scale=1.7,trim = 0mm 3mm 22mm 0mm, clip]{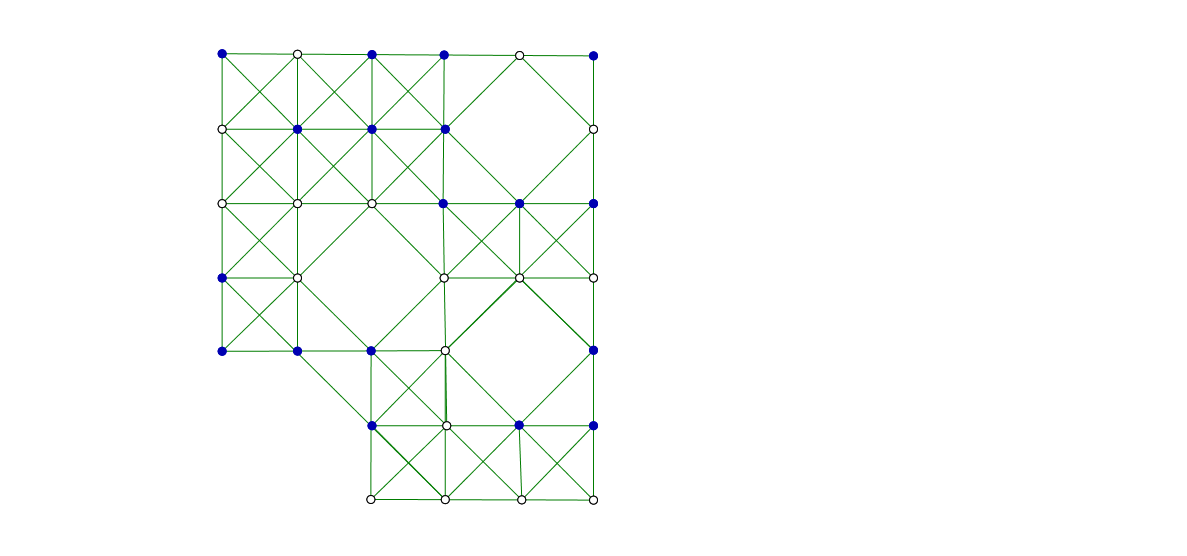}
\end{center}
\caption{The graph corresponding to the matrix of the example. Blue nodes define 1's in the matrix and white nodes 0's. The nodes corresponding to entries other than 0 and 1 have been removed, since they cannot be visited.}
\label{Figure 2}
\end{figure}

It is a well known result that the Hamiltonian Path problem in node-induced subgraphs of the grid is a NP-complete problem\cite{itai}. We will, now, show that the Longest Path problem in the grid is a NP-complete problem.

\subsection{The Longest Path in a Grid problem is NP-complete}

First of all, the Longest Path in a Grid problem belongs to NP, since we can verify in polynomial time that when given a node-induced subgraph of a grid with unit lengths on all edges, a positive integer $k\leq|V|$ and a set of 2 distinct nodes of the grid, we can verify in linear time that a path is simple having length not less than k and connecting the 2 given vertices to each other. Also, using a 'black box' solving this problem, we would be able to reply for the Hamiltonian path in a grid problem which is NP-complete as described before, provided that we were given that $k=|V|$, so we can easily infer that this problem is NP-complete. 

Finally, it is easy to induce that the problem whether there is a simple path of length not less than k for any pair of nodes in a grid, is, also, NP-complete, since, we are actually applying the Longest Path in a Grid problem $\mathcal{O}(n^{2})$ times.   

\subsection{The 0-1 String to Matrix problem is NP-complete}

It is obvious, that the 0-1 String to Matrix problem belongs to NP, since we can easily verify (in linear time) that a path in the matrix is identical to a provided sequence, visiting 0's and 1's no more than once each.

Now, we need a way to represent the matrix as a graph. Since all available moves can be made horizontally, vertically, or diagonally and only to adjacent entries, we  may represent the matrix as a diagonal grid, where for 0's we place white nodes and for 1's blue nodes, as in figure 1, representing the matrix of the previous example. For entries other than 0 and 1, we remove the corresponding nodes, obtaining a node-induced subgraph of the diagonal grid, since these entries cannot be visited.

Since rectangular grids are bipartite graphs, they may be colored using only 2 colors, where adjacent vertices share different colors. So, we may color a random node-induced subgraph of a grid with 2 colors, blue and white, starting from the upper left entry and continuing in the same row, assigning different colors to adjacent entries and afterwards to the second line in reverse order and so on, until we reach the last row. Each time we pass through a vacant area (where a set of nodes was removed, obtaining the node-induced rectangular grid) we continue assigning colors to the rest of the graph nodes, as that area did not exist. The crucial observation is that all nodes observed diagonally (even though not connected directly of course) share the same colors.

We will now reduce the Longest Path in a Grid problem (for all pair of nodes) to the 0-1 String to Matrix problem. If we were able to solve the 0-1 String to Matrix problem for a matrix whose each column or row would not have two or more 0's or 1's in two consecutive positions (in the exact same column or row respectively) and for the specific sequence of non-repeating identical digits in two positions in a row 010101...01, then we would be able to solve the Longest Path in a Grid problem (for all pairs of nodes) in the corresponding node-induced rectangular grid. We note that we choose a sequence of non-repeating identical digits according to how many blue and white nodes we, finally, have in the coloring we have assigned, as described before. So, if we have p blue nodes and p+q white ones, we choose the sequence of 2p digits 10101010...10, taking under account only the least of the two numbers of digits. The key point is that, because  of the sequence 101010...10, the diagonal edges of the graph representing the matrix, cannot be traversed since they share the same color, assuming we color the graph as described before; in other words we obtain a node-induced grid by removing them, where k equals the length of the sequence. The NP-completeness of the 0-1 String to Matrix problem follows.

\begin{figure}
\begin{center}
\includegraphics[scale=1.7,trim = 12mm 3mm 5mm 2mm, clip]{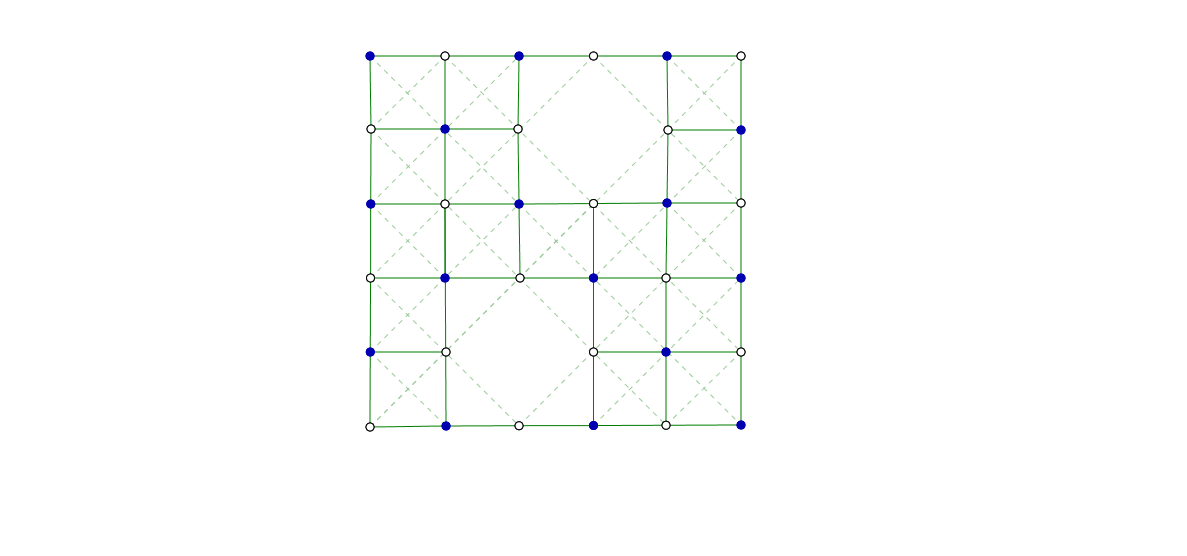}
\end{center}
\caption{The graph corresponding to a 6x6 matrix where the entries correspond to 1's for the blue nodes and to 0's for the white ones and where we, also, have 2 entries with digits other than 0 and 1. We have 16 '1' and 18 '0' entries, so we choose the 32-digit sequence 10101010101010101010101010101010, making the diagonals useless for a possible simple path matching this sequence.}
\label{Figure 3}
\end{figure}

Furthermore, we can restrict our problem to only square matrices, by padding the matrix with digits other than 0 and 1, obtaining a square lxl matrix, where l=max(n,m).

\section{Open Problems}
The restricted case where no digits other than 0 and 1 are allowed in the square nxn matrix, is an interesting special open case of this problem, assuming the sequence has length equal to $n^{2}$. Also, an other interesting problem would be the more general case of that, whether it is possible, after breaking the sequence in maximum substrings of identical digits, to form a matching path in the matrix (formed by only 0's and 1's) at any order of the substrings.




\bibliographystyle{plain}	
\bibliography{myrefs}

\begin{thebibliography}{1}

\bibitem{eppstein}
http://www.ics.uci.edu/$\sim$eppstein/200-f01.pdf.

\bibitem{garey}
Michael~R. Garey and David~S. Johnson.
\newblock {\em Computers and Intractability. A guide to the Theory of
  NP-Completeness}.
\newblock W.H. Freeman and Company, 1979.

\bibitem{itai}
Alon Itai, Christos~H. Papadimitriou, and Jayme~Luiz Szwarcfiter.
\newblock Hamilton paths in grid graphs.
\newblock {\em SIAM J. Comput.}, 11(4):676--686, 1982.

\bibitem{pap}
Christos Papadimitriou.
\newblock {\em Computational Complexity}.
\newblock Addison-Wesley, 1994.

\end{thebibliography}
\end{document}